\documentclass[12pt]{article}
\usepackage{graphicx}
\usepackage{amsmath,amssymb}
\newtheorem{theorem}{Theorem}[section]

\newtheorem{proposition}[theorem]{Proposition}

\begin{document}



\title{An improved method for the computation of the Moore-Penrose inverse matrix}

\author{Vasilios N. Katsikis\\General Department of Mathematics\\Technological Education Institute of Piraeus,\\Aigaleo, 12244 Athens,
Greece, vaskats@gmail.com\\ \\Dimitrios Pappas\\Athanassios Petralias\\ Department of Statistics\\
Athens University of Economics and Business\\ 76 Patission Str,
10434, Athens, Greece\\pappdimitris@gmail.com, petralia@aueb.gr}
\maketitle

\begin{abstract}
In this article we provide a fast computational method in order to
calculate the Moore-Penrose inverse of singular square matrices and
of rectangular matrices. The proposed method proves to be much faster and has
significantly better accuracy than the already proposed methods, while works for full
and sparse matrices.
\end{abstract}

\textit{Keywords}: Moore-Penrose inverse matrix, tensor-product
matrix.



\section{Introduction} \label{intro-sec}
Let $T$ be a  $n\times n$ real matrix. It is known that when $T$ is
singular, then its unique generalized inverse $T^{\dagger}$ (known
as the Moore- Penrose inverse) is defined. In the case when $T$ is a
real $m\times n$ matrix, Penrose showed that there is a unique
matrix satisfying the four Penrose equations, called the generalized
inverse of $T$. A lot of work concerning generalized inverses has
been carried out, in finite and infinite dimension (e.g.,
\cite{Israel,RAKH}).

In a recent article \cite{Kats}, the first two authors provided a
new method for the fast computation of the generalized inverse of
full rank rectangular matrices and of square matrices with at least
one zero row or column. In order to reach this goal, a special type
of tensor product of two vectors was used, usually defined in
infinite dimensional Hilbert spaces. In this work we extend our
method so that it can be used in any kind of matrix, square or
rectangular, full rank or not. The numerical experiments show that
the proposed method is competitive in terms of accuracy and much
faster than the commonly used methods, and can be also used for
large sparse matrices. \\There are several methods for computing the
Moore-Penrose inverse matrix (cf. \cite{Israel}). One of the most
commonly used methods is the Singular Value Decomposition (SVD)
method. This method is very accurate but also time-intensive since
it requires a large amount of computational resources, especially in
the case of large matrices. On a recent work, Toutounian and Ataei
\cite{Tout} presented an algorithm based on the conjugate
Gram-Schmidt process and the Moore-Penrose inverse of partitioned
matrices, the  \verb CGS-MPi   algorithm and they conclude that this
algorithm is a robust and efficient tool for computing the
Moore-Penrose inverse of large sparse and rank deficient matrices.

In the present manuscript, we construct a very fast and reliable
method (see the \verb qrginv  function in the Appendix) in order to
estimate the Moore-Penrose inverse matrix. The computational effort
required for the \verb qrginv  function (see Figure 1 and Tables 1,
4) in order to obtain the generalized inverse is substantially
lower, particularly for large matrices, compared to those provided
by the SVD and the method presented in \cite{Tout}
   (\verb CGS-MPi   algorithm).  In addition, we obtain reliable and very accurate
approximations in each one of the tested cases (Tables 2, 3 and 4).
In order to test this algorithm, we have used random singular
matrices (see subsection 4.1) as well as a collection of singular
test matrices (see subsection 4.2) with ``large" condition number
(ill-conditioned matrices) from the Matrix Computation Toolbox (see
\cite{HIGH}). We also tested the proposed method on sparse matrices
from the Matrix Market collection \cite{MAT} (see subsection 4.3).
In what follows, we make use of the high-level language Matlab both
for calculations of the generalized inverse matrix, as well as for
testing the reliability of the obtained results. Specifically, the
Matlab 7.4 (R2007a) Service Pack 3 version of the software was used
on an Intel Core i7 920 Processor system running at 2.67GHz with 6
GB of RAM memory using the Windows XP Professional 64-bit Operating
System.

\section{Preliminaries and notation}
We shall denote by $\mathbb{R}^{m\times n}$ the linear space of all
$m\times n$ real matrices. For $T \in \mathbb{R}^{m\times n}$, the
generalized inverse $T^{\dagger}\in \mathbb{R}^{n\times m}$ (known
as the Moore- Penrose inverse) is the unique matrix that satisfies
the following four Penrose equations:
$$TT^{\dagger}=(TT^{\dagger})^*,\qquad T^{\dagger}T=(T^{\dagger}T)^*,\qquad
   TT^{\dagger}T=T,\qquad T^{\dagger}TT^{\dagger}=T^{\dagger},$$ where $T^{*}$ denotes the transpose matrix of $T.$
The number $r = \dim \mathcal{R}(T)$ is called the rank of $T$ and
$\langle\: ,\rangle$ denotes the usual inner-product in
$\mathbb{R}^n.$

According to \cite{Ring}, for each $x\in \mathbb{R}^k,$ we consider
the mapping
$$e\otimes f:\mathbb{R}^k\rightarrow \mathbb{R}^k \textit{ with }
(e\otimes f)(x)=\langle x,e\rangle f.$$ Assume that
$\{e_1,...,e_n\}$ and $\{f_1,...,f_n\}$ are two collections of
orthonormal vectors and linearly independent vectors of
$\mathbb{R}^k$ with $n<k$, respectively. Then, every rank-$n$
operator $T$ can be written in the  form $T=\sum_{i=1}^n e_i\otimes
f_i.$ We shall refer to this type of tensor product as the
\textit{tensor-product of the collections} $\{e_1,...,e_n\}$ and
$\{f_1,...,f_n\}.$ The adjoint operator $T^*$ of $T$ is the rank-$n$
operator $T^*=\sum_{i=1}^n f_i\otimes e_i$.

The tensor-product of two collections of vectors, as defined above,
is a linear operator, therefore, it has a corresponding matrix
representation $T$. We shall refer to this matrix $T$ as the
\textit{tensor-product matrix} of the given collections. In order to
compute  the Moore-Penrose inverse of the corresponding
tensor-product matrix, we use the following theorem,
\begin{theorem}[\cite{kar1}, Theorem 3.2]\label{T32}
Let $\mathcal{H}$ be a Hilbert space. If $T=\sum_{i=1}^n e_i \otimes
f_i$ is a rank-$n$ operator then its generalized inverse is also a
rank-$n$ operator and for each $x\in \mathcal{H}$, it is defined by
the relation
$$T^{\dagger}x=\sum_{i=1}^n \lambda_i (x)e_i,$$ where the functions $\lambda_i
$ are the solution of an appropriately defined $n\times n$ linear
system.
\end{theorem}
\section{The computational method}
In \cite{Kats}, based on theorem 1, the authors developed an
algorithm (the \verb ginv  function) for computing the generalized
inverse of full rank matrices, and of square matrices with at least
one zero row or column and the rest of the matrix full rank. In
other words, our main concern was to calculate the corresponding
$\lambda_i$ in the expansion
$$T^{\dagger}x=\sum_{i=1}^n \lambda_i (x)e_i,$$ where $\{e_1,...,e_n\}$ are the
first $n$ vectors of the standard basis of $\mathbb{R}^k$, in order
to provide the generalized inverse $T^{\dagger}$.

To extend this result for any kind of matrix, we will make use of
the QR factorization, as well as the reverse order law for
generalized inverses. The following proposition is a restatement of
a part of R. Bouldin's theorem \cite{boul} which holds for operators
and matrices (see also \cite{Grev}, \cite{Izum}).
\begin{proposition}\label{bouldin}
Let $A,B$ be bounded operators on $\mathcal{H}$ with closed range.
Then $(AB)^\dagger=B^\dagger A^\dagger $ if and only if the
following three conditions hold:

i) The range of $AB$ is closed,

ii) $A^\dagger A$ commutes with $BB^*$,

iii) $BB^\dagger$ commutes with $A^*A$.
\end{proposition}
Using the above proposition, we have the following result, which can
be found also, in a similar form but with a different proof, in
\cite{Israel}.
\begin{proposition}\label{prop} Let $A=QR$ be the QR factorization of $A$.
Then, $A^\dagger = R^\dagger Q^*$.
\end{proposition}
Proof:  We must prove that the conditions of Bouldin's theorem hold.
The first condition always holds, since in the case of matrices the
range is always closed. For the second condition, it is easy to see
that since  $Q$ is a unitary matrix, $Q^\dag = Q^*=Q^{-1}$ and so \[
Q^\dagger QRR^* =Q^{-1}QRR^*= IRR^*=RR^*I=RR^*Q^\dagger Q.\] The
third condition can be proved in a similar way.
\\
Using the QR factorization, the matrix $R$ is upper triangular but
not necessarily of full rank. So a variant of the QR method must be
used, the QR with column pivoting as described in the following form
from \cite{wat}:
 \begin{theorem}[\cite{wat}, Theorem 3.3.11]
Let $A\in \mathbb{R}^{n\times m}$ matrix, with $ rank(A)=r >0$. Then
there exist matrices $\hat{A}, Q, R$ such that $\hat{A}$ is obtained
by $A$ by permuting its columns, $Q \in \mathbb{R}^{n\times n}$ is
orthogonal, $R = \left[ {\begin{array}{*{20}c}
   R_{11} & R_{12}  \\
   0 & 0  \\
\end{array}} \right] \in\mathbb{R}^{ n\times m}$, $R_{11} \in \mathbb{R}^{r\times r}$ is nonsingular and
upper triangular and $\hat{A}= Q R$.
\end{theorem}

Using the above theorem, we have that $AP =QR$, where $P$ is a
permutation matrix (therefore unitary). By proposition \ref{prop} we
have that $A^\dagger=P R^\dagger Q^*$.

To calculate the rank  of $R_{11}$, one needs only the number of its
columns that have at least one value above a tolerance level in
absolute terms. This tolerance is set equal to $10^{-5}$, which is
also used by Toutounian and Ataei \cite{Tout}, and turns out to
provide accurate results. The implementation of all the above ideas
are presented in the \verb qrginv  function (see the Appendix).

\section{Numerical experiments}
In this section we perform numerical experiments comparing Matlab's
 \verb pinv  function, Toutounian and Ataei' s method
\cite{Tout} ( \verb CGS-MPi  algorithm) and the proposed method
\verb qrginv  function. Testing of  \verb pinv , \verb CGS-MPi  and
\verb qrginv  was performed separately for random singular matrices
and for singular test matrices with ``large" condition number from
the Matrix Computation Toolbox (see \cite{HIGH}). We also tested the
proposed method in sparse matrices and we obtained  very fast and
accurate results.

\subsection{Random singular matrices}
We are comparing the performance of the proposed method (\verb
qrginv)
  to that of the SVD method used by Matlab (\verb pinv ) and the method proposed in \cite{Tout} (\verb CGS-MPi  algorithm),
  on a series of random singular matrices with rank  $2^{n},\,\,n=8,9,10,11,12$. In addition, the
accuracy of the results is examined with the matrix 2-norm in error
matrices corresponding to the four properties characterizing the
Moore-Penrose inverse.

\begin{table}[h]
\label{v1} \caption{Computational Time Results; Random singular matrices}
\scriptsize{\begin{tabular}{llll}
Rank& pinv & CGS-MPi & qrginv \\
  \hline
  \\
$2^{8}$& 0.496& 1.281 &0.0247\\
$2^{9}$& 1.185 & 15.183 & 0.104 \\
$2^{10}$ & 14.066 & 448.172 & 1.556 \\
$2^{11}$ & 152.384 & 11374.47 & 9.259 \\
$2^{12}$ & 1081.81 & - & 45.061\\
\\
\hline
\multicolumn{4}{l}{Notes: Time is measured in seconds. \texttt{CGS-MPi} }\\
\multicolumn{4}{l}{algorithm was not able to produce numerical}\\
\multicolumn{4}{l}{results for a rank $2^{12}$ matrix, even after 3 days.}
\end{tabular}}
\end{table}

In Table 1 we can see the time efficiency for the same matrices of the proposed method, the CGS-MPi algorithm and the \verb pinv  command used by Matlab. The \verb qrginv  method needs about 4\% up to 11\% the corresponding time needed by the \verb pinv  function to calculate the Moore-Penrose inverse, depending on the matrix. On the other hand the \verb CGS-MPi  turns to be computantionally demanding requiring from 50 times up to more than 1200 times the corresponding time needed by \verb qrginv . Furthermore, the larger the rank of the matrix, the greater this difference, so that for a matrix with rank $2^{12}$, the \verb CGS-MPi  algorithm was not able to produce a numerical result, even after 3 days running, while \verb qrginv  needs only up to 45 seconds.

\begin{table}[h]
\label{v1} \caption{Error Results; Random singular matrices}
\scriptsize{\begin{tabular}{lcllll}
Method&Rank&\scriptsize{$\|TT^{\dagger}T-T\|_2$} &\scriptsize{$\|T^{\dagger}TT^{\dagger}-T^{\dagger}\|_2$}&\scriptsize{$\|TT^{\dagger}-(TT^{\dagger})^*\|_2$}&\scriptsize{$\|T^{\dagger}T-(T^{\dagger}T)^*\|_2$} \\
  \hline
  \\
   \scriptsize{pinv}&& $1.21  \times 10^{-12}$ &$5.58 \times 10^{-13}$&$2.98 \times 10^{-13}$&$3.46 \times 10^{-13}$\\
 \scriptsize{CGS-MPi}&$2^{8}$& $2.13  \times 10^{-8}$ &$7.44 \times 10^{-7}$&$2.01 \times 10^{-11}$&$5.29\times 10^{-7}$\\
     \scriptsize{qrginv}&&$  1.44 \times 10^{-13}$&0&0 &0\\  & & & & &\\
\\

\scriptsize{pinv}&& $3.07 \times 10^{-12}$ &$2.31 \times 10^{-12}$&$1.49 \times 10^{-12}$&$7.33 \times 10^{-13}$\\
\scriptsize{CGS-MPi}&$2^{9}$& $1.62  \times 10^{-8}$ &$9.79 \times 10^{-9}$&$9.68 \times 10^{-11}$&$5.18 \times 10^{-9}$\\
     \scriptsize{qrginv}&&$  3.62 \times 10^{-13}$&0&0 &0\\  & & & & &\\
\\

\scriptsize{pinv}&& $8.24  \times 10^{-12}$ &$3.18 \times 10^{-12}$&$1.63 \times 10^{-12}$&$1.62 \times 10^{-12}$\\
\scriptsize{CGS-MPi}&$2^{10}$& $1.57  \times 10^{-8}$ &$3.70 \times 10^{-7}$&$2.90 \times 10^{-10}$&$1.65 \times 10^{-7}$\\
     \scriptsize{qrginv}&&$  7.47 \times 10^{-13}$&0&0 &0\\  & & & & &\\
\\

\scriptsize{pinv}&& $4.11  \times 10^{-11}$ &$2.90 \times 10^{-11}$&$9.32 \times 10^{-12}$&$9.45 \times 10^{-12}$\\
\scriptsize{CGS-MPi}&$2^{11}$& $6.99  \times 10^{-7}$ &$4.19 \times 10^{-5}$&$5.66 \times 10^{-9}$&$1.15 \times 10^{-5}$\\
     \scriptsize{qrginv}&&$  6.31 \times 10^{-12}$&0&0 &0\\  & & & & &\\
\\

\scriptsize{pinv}&& $1.19  \times 10^{-10}$ &$7.52 \times 10^{-10}$&$2.43 \times 10^{-11}$&$5.69 \times 10^{-11}$\\
\scriptsize{CGS-MPi}&$2^{12}$& - &-&-&-\\
     \scriptsize{qrginv}&&$  1.01 \times 10^{-11}$&0&0 &0\\  & & & & &\\
\\
\hline
\multicolumn{6}{l}{Notes: Zeros (0) denote numbers close to zero under 64-bit IEEE double
precision arithmetic. The \texttt{CGS-MPi} }\\
\multicolumn{6}{l}{algorithm was not able to produce numerical results for matrix with rank $2^{12}$ even after 3 days running.}
\end{tabular}}
\end{table}

In Table 2, the accuracy of the proposed method is tested based on the 2-norm errors. It is evident that the proposed method (\verb qrginv)  produced a reliable approximation in all the tests that were conducted. The associated errors are greatly lower than the corresponding of the other methods, while in certain (many) cases are equal to zero, under the 64-bit IEEE double precision arithmetic in Matlab. Therefore, the proposed method allows us for a both fast and accurate computation of the Moore-Penrose inverse matrix.

\subsection{Singular test matrices}
In this section we use a set of singular test matrices
that includes 13 singular matrices, of size $200\times 200$, obtained from
the function
 \verb matrix   in the Matrix Computation Toolbox \cite{HIGH}(which
 includes test matrices from Matlab itself). The condition number of
 these matrices range from order $10^{15}$ to $10^{135}$. For comparative purpose
 we also apply as in the previous section, Matlab's \verb pinv  function which implements
 the SVD method and \verb CGS-MPi  algorithm.  Since these matrices are of relatively small size and so as to measure the time needed
 for each algorithm to compute the Moore-Penrose inverse accurately, each
 algorithm runs 100 distinct times. The reported time is the mean
 time over these 100 replications. The error results are presented in Table 3,
 while the time responses are shown in Figure 1.

 The results show a clear ordering of the three methods for this set of test
 problems, with \verb qrginv  in the first place, followed by
\verb pinv   and then by \verb CGS-MPi  algorithm. The worst cases
for \verb pinv,  with comparatively large errors for the
$\|T^{\dagger}TT^{\dagger}-T^{\dagger}\|_2$ norm, are the lotkin,
prolate, hilb, and vand matrices. The \verb CGS-MPi  algorithm has
very large errors in the cases of Kahan, lotkin, prolate, hilb,
magic and vand matrices, for all tested norms. Moreover, the
algorithm failed to produce a numerical result for the chow matrix,
since the computed Moore-Penrose inverse included only NaNs.
Nevertheless, we observe that there are also cases (in certain
norms) that the \verb CGS-MPi  algorithm has smaller error than
\verb pinv. On the other hand, the proposed \verb qrginv  method
gives very accurate results for all matrices and proves to be
overally more efficient than the other two methods.

\begin{center}
\begin{figure}[h!] \label{f2}
\includegraphics[scale=8, angle=0,width=1\textwidth, height=2.6in]{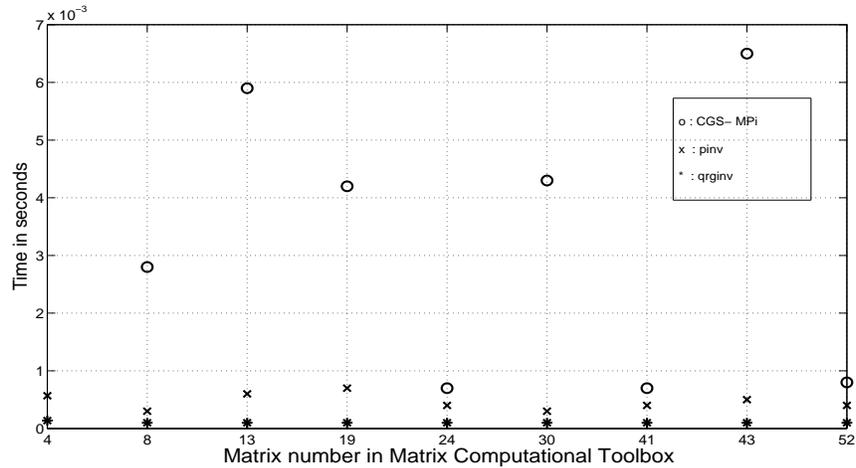}\caption{
Time efficiency for the\texttt{ pinv}, \texttt{qrginv} functions and
the \texttt{CGS-MPi} algorithm }
\end{figure}
\end{center}

\begin{table}[h!]
\label{v1} \caption{Error Results; Singular test matrices}
\scriptsize{\begin{tabular}{lllllll}
\tiny{Matrix}&\tiny{Condition number}&\tiny{Method}&\tiny{$\|TT^{\dagger}T-T\|_2$} &\tiny{$\|T^{\dagger}TT^{\dagger}-T^{\dagger}\|_2$}&\tiny{$\|TT^{\dagger}-(TT^{\dagger})^*\|_2$}&\tiny{$\|T^{\dagger}T-(T^{\dagger}T)^*\|_2$} \\
  \hline
  \\

\tiny{chow}& \tiny{$8.01849\times 10^{135}$} &\tiny{pinv}&
$4.741\times 10^{-13}$  &$1.896\times 10^{-13}$ &$2.313\times
10^{-13}$  &$2.261\times 10^{-13}$ \\ \tiny{(r=199)} & &\tiny{CGS-MPi} & -&-&-&-\\
& &\tiny{qrginv}&$1.691\times 10^{-13}$&0&0&0\\
 \\

\tiny{cycol}&\tiny{$2.05\times 10^{48}$ }&\tiny{pinv}& $1.539  \times 10^{-13}$ &$1.637 \times 10^{-16}$&$4.285 \times 10^{-15}$&$3.996 \times 10^{-15}$\\ \tiny{(r=50)}
& &\tiny{CGS-MPi} & $3.507  \times 10^{-14}$ &$6.365
\times10^{-17}$&$6.691 \times 10^{-16}$&$1.236 \times 10^{-15}$ \\
& &\tiny{qrginv}&$  1.262 \times 10^{-15}$&0&0 & 0\\
\\

\tiny{gearmat}& \tiny{$3.504\times 10^{17}$}
&\tiny{pinv}&$1.899\times 10^{-14}$&$3.033\times
10^{-13}$&$8.063\times 10^{-14}$ &$7.561\times 10^{-14}$\\ \tiny{(r = 199)}
& &\tiny{CGS-MPi} & $2.379 \times 10^{-13}$ &$3.392
\times10^{-13}$&$6.113 \times 10^{-14}$&$5.478 \times 10^{-13}$ \\
& & \tiny{qrginv}&$  1.923 \times 10^{-14}$&0&0 & 0\\
 \\

\tiny{kahan}& \tiny{$2.30018\times 10^{24}$}
&\tiny{pinv}& $1.432  \times 10^{-13}$ &0&0&$9.070 \times 10^{-14}$\\ \tiny{(r = 199)}
& &\tiny{CGS-MPi} & $4.011 \times 10^{+18}$ &$4.921
\times10^{+35}$&$262.069 $&$5.765 \times 10^{+18}$ \\
& &  \tiny{qrginv}&$  6.964 \times 10^{-15}$&0&0 & 0\\
 \\

\tiny{lotkin}&  \tiny{$8.97733\times 10^{21}$}
&\tiny{pinv}& $0.0001$ &$103.5\times 10^{+6}$&$0.001$&$0.0011$\\ \tiny{(r = 19)}
& &\tiny{CGS-MPi} & $29.59 \times 10^{+7}$ &$1.252
\times10^{+14}$&$0$&$74.5 \times 10^{+6}$ \\
& &  \tiny{qrginv}&$  3.546 \times 10^{-11}$&0&0 & 0\\
 \\

\tiny{prolate}&  \tiny{$5.61627\times 10^{17}$}
&\tiny{pinv}& $0.0002$ &$7.398\times 10^{9}$&$0.0096$&$0.0081$\\ \tiny{(r = 117)}
& &\tiny{CGS-MPi} & $2.218 \times 10^{+15}$ &$7.885\times10^{+30}$&$0.625$&$4.796 \times 10^{+15}$ \\
& &\tiny{qrginv}&$  6.588 \times 10^{-11}$&0&0 & 0\\
 \\

\tiny{hilb}& \tiny{$1.17164\times 10^{19}$} &\tiny{pinv}& $0.0001$ &$3.184\times 10^{9}$&$0.0033$&$0.0051$\\ \tiny{(r = 20)}
& &\tiny{CGS-MPi} & $39.6 \times 10^{+5}$ &$2.811
\times10^{+12}$&$5.533 \times10^{-10}$&$30.95 \times 10^{+5}$ \\
& &\tiny{qrginv}&$  6.941 \times 10^{-12}$&0&0 & 0\\
 \\

\tiny{magic}& \tiny{$4.92358\times 10^{19}$} &\tiny{pinv}&$0$ &$2.491\times 10^{-19}$&$9.701\times 10^{-14}$&$4.236\times 10^{-14}$\\ \tiny{(r = 3)}
& &\tiny{CGS-MPi} & $1.755 \times 10^{+18}$ &$8.333
\times10^{+19}$&$293.29$&$6.329 \times 10^{+14}$\\
& &\tiny{qrginv}&$  4.566 \times 10^{-14}$&0&0 & 0\\
 \\

\tiny{vand}&\tiny{$1.16262\times 10^{20}$}&\tiny{pinv}&$0.0003$ &$179.8\times 10^{+6}$&$0.0022$&$0.0019$\\ \tiny{(r = 34)}
& &\tiny{CGS-MPi} & $6.498 \times 10^{+11}$ &$5.299
\times10^{+20}$&0&$2.39 \times 10^{+11}$\\
 & &\tiny{qrginv}&$  5.553 \times 10^{-11}$&$0$&$0$ & 0\\
 \\

\hline
\end{tabular}}
Notes: Zeros (0) denote numbers close to zero under 64-bit IEEE double
precision arithmetic. In parenthesis is denoted the rank (r) of each matrix. The \verb CGS-MPi  algorithm was not able to produce numerical results for the chow matrix.
\end{table}

\subsection{Matrix Market sparse matrices}
In this section we test the proposed algorithm on sparse matrices,
 from the Matrix Market collection \cite{MAT}. We follow the same
 method and the same matrices as in \cite{Tout}, while the matrices are taken as rank deficient: $A\_Z = [A\ \  Z]$, where A is one of the chosen matrices, Z is a zero matrix of order $m \times 100$ and $m$ takes values
shown in Table 4, as in \cite{Tout}.

The proposed algorithm is tested against the \verb CGS-MPi  algorithm, since this method is proposed by Toutounian and Ataei
\cite{Tout} as suitable for large and sparse matrices. The \verb pinv  function of Matlab is not applicable in sparse matrices and thus ommited. We observe that in sparse matrices as well, the proposed method seems to greatly outperform the \verb CGS-MPi  algorithm, both in terms of speed and accuracy.

\begin{table}[h!]
\label{v1} \caption{Error and Computational Time Results; Matrix Market sparse matrices}
\scriptsize{\begin{tabular}{lllllll}
\tiny{Matrix}&\tiny{Time}&\tiny{Method}&\tiny{$\|TT^{\dagger}T-T\|_2$} &\tiny{$\|T^{\dagger}TT^{\dagger}-T^{\dagger}\|_2$}&\tiny{$\|TT^{\dagger}-(TT^{\dagger})^*\|_2$}&\tiny{$\|T^{\dagger}T-(T^{\dagger}T)^*\|_2$} \\
  \hline
  \\
\tiny{WELL1033\_Z}& $0.0716$ &\tiny{CGS-MPi}& $2.801\times
10^{-11}$ &0 &$9.307\times 10^{-12}$  &$1.303\times 10^{-10}$
\\\tiny{($m =1033$)}& $0.0111$&\tiny{qrginv} & $1.142  \times
10^{-12}$ &0&0&0 \\
\\

\tiny{WELL1850\_Z}&$0.4626$  &\tiny{CGS-MPi}& $2.184\times
10^{-11}$ &$1.136\times 10^{-10}$ &$7.135\times 10^{-12}$
&$7.236\times 10^{-11}$
\\\tiny{($m =1850$)}&$0.0315$ &\tiny{qrginv} & $1.89  \times 10^{-12}$ &0&0&0 \\
\\

\tiny{ILCC1850\_Z}& $0.4546$ &\tiny{CGS-MPi}& $61.106\times
10^{+4}$ &$1.687\times 10^{+9}$ &$2.637\times 10^{-10}$
&$50.971\times 10^{+4}$
\\\tiny{($m =1850$)}&$0.0315$ &\tiny{qrginv} & $5.809  \times 10^{-11}$ &0&0 & 0 \\
  \\

\tiny{GR-30-30\_Z}& $0.9618$ &\tiny{CGS-MPi}& $4.765\times
10^{-10}$ &$5.91\times 10^{-11}$ &$2.206\times 10^{-11}$
&$6.746\times 10^{-10}$
\\\tiny{($m =900$)}&$0.0336$ &\tiny{qrginv} & $7.321  \times 10^{-12}$ &0&0 & 0 \\
  \\

\tiny{WATT1\_Z}& $1.2687$ &\tiny{CGS-MPi}& 8 &$101.6\times 10^{+6}$
&$7.743\times 10^{-16}$  &$27.783$ \\\tiny{($m =1856$)}&$0.0064$
&\tiny{qrginv} & 0 &0&0 & 0 \\
\\
\hline
\end{tabular}}
 Notes: Zeros (0) denote numbers close to zero under 64-bit IEEE double
precision arithmetic. In parenthesis is denoted the row size of each matrix ($m\times 100$). Time is measured in seconds.
\end{table}

\section{Concluding Remarks}
We proposed a new method for calculating the Moore-Penrose inverse
of singular square, rectangular, full or sparse matrices. It is
apparent that the proposed method provides a substantially faster
numerical way for calculating the Moore-Penrose inverse of a given
matrix (see Figure 1 and Tables 1, 4). Also, it is evident, from
Tables 2, 3 and 4 that the proposed function
 (\verb qrginv )provides a far more reliable approximation in all the tests,
compared to other existing methods.

\appendix
\begin{flushleft}
\textbf{Appendix: Matlab code of the qrginv function}
\end{flushleft}
\begin{verbatim}
function qrginv = qrginv(B) [N,M] = size(B); [Q,R,P] = qr(B);
r=sum(any(abs(R)>1e-5,2)); R1 = R(1:r,:); R2 = ginv(R1); R3 = [R2
zeros(M,N-r)]; A = P*R3*Q'; qrginv = A;
\end{verbatim}

In case the matrix of interest is sparse, the third line of the code
is replaced with
\begin{verbatim}
[Q,R,P] = spqr(B);
\end{verbatim}
since the function \texttt{qr} embeded in Matlab does not support
sparse matrices under this format. The function \texttt{spqr} is
part of the SuiteSparse toolbox, built by Professor Timothy A.
Davis, University of Florida and can be downloaded electronically
from \emph{http://www.cise.ufl.edu/research/sparse/SuiteSparse/}.



\end{document}